\documentclass[12pt]{article}
\title{On integers $n$ for which $\sigma(2n+1)\geq \sigma(2n)$}
\author{Mits Kobayashi\\ Dartmouth College\\ mitsuo.kobayashi@dartmouth.edu\and  Tim Trudgian\footnote{Supported by Australian Research Council Future Fellowship FT160100094.} \\
School of Science\\ University of New South Wales Canberra, Australia \\
  t.trudgian@adfa.edu.au}
\usepackage{url}
\usepackage{enumerate}
\usepackage{amsthm}
\usepackage{amsmath}
\usepackage{comment}
\usepackage{fullpage}
\usepackage{refcheck}
\usepackage{amssymb}
\usepackage{booktabs}
\usepackage{relsize}
\usepackage{mathrsfs}

\newtheorem{thm}{Theorem}

\newtheorem{lem}{Lemma}

\newcommand{\bbZ}{\mathbb Z}
\newcommand{\bbN}{\mathbb N}

\DeclareMathOperator \dens {\mathbf{d}}

\newcommand{\upperbound}{0.0549445}
\newcommand{\lowerbound}{0.0539171}
\newcommand{\shortubd}{0.055}
\newcommand{\shortlbd}{0.053}
\newcommand{\deleglise}{Del\'eglise}

\begin{document}
\maketitle
\begin{abstract}
\noindent
We show that the natural density of positive integers $n$ for which $\sigma(2n+1)\geq \sigma(2n)$ is between $\shortlbd$ and $\shortubd$.
\end{abstract}

\section{Introduction}
Let $\sigma(n)$ denote the sum of divisors function. While its average value is well-behaved (see, e.g.\ \cite[\S 18.3]{HW}), the local behavior of $\sigma(n)$ is, as with many interesting arithmetical functions, erratic. Consider, for example, a result from 
Erd\H{o}s, Gy\H{o}ry, and Papp \cite{EGP} (see also \cite[p.\ 89]{Sandor}) that says that the chain of inequalities
$$ \sigma(n+ m_{1}) > \sigma(n+ m_{2}) > \sigma(n+ m_{3}) > \sigma(n+ m_{4})$$
holds for infinitely many $n$, where the $m_{i}$ are any permutations of the numbers $1,2,3,4$.

We consider here the problem of counting those $n$ such that $\sigma(2n+1) \geq \sigma(2n)$. When $2n+1$ is prime the left side is $2n+2$ whereas the right side is at least $2n+1 + n + 2 = 3n+3$. This shows that the inequality is false infinitely often. Empirically, it appears to be false very frequently.
Let $B$ be the set of natural numbers $n$ such that $\sigma(2n+1)\geq \sigma(2n)$ and let $B(x)$ be the number of those $n$ in $B$ with $n\leq x$. From Table \ref{strauss} one may be tempted to conjecture that $B(x)/x \sim 0.0546\ldots$ .

\begin{table}[ht]
\caption{Proportion of integers $n\leq x$ with $\sigma(2n+1)\geq \sigma(2n)$}
\centering
\begin{tabular}
{c c}
\hline\hline
$x$& Proportion\\[0.5ex]\hline
$10^3$& $0.06$\\
$10^4$& $0.0551$\\
$10^5$& $0.549$\\
$10^6$& $0.054603$\\
$10^7$& $0.0546879$\\
$10^8$& $0.0546537$\\
$10^{9}$ & $0.054665173$\\
    \hline\hline
  \end{tabular}
\label{strauss}
\end{table}

Laub \cite{Laub} posed the question of estimating the size of $B(x)/x$. Mattics \cite{Mattics} answered this, and records a remark of Hildebrand that $\lim_{x\rightarrow\infty}B(x)/x$ exists. We will call this limit the natural density of $B$, denoted $\dens B$. Although Mattics was not able to calculate this density, he was able to establish the existence of constants $\lambda$ and $\mu$ with $0<\lambda<\mu<1$ such that $\lambda x<B(x)<\mu x$ for $x$ sufficiently large. Specifically, he showed that one could take $\lambda =1/3000$ and $\mu= 25/28$.

We refine Mattics' result and prove the following.
\begin{thm}\label{flower}
Let $B= \{n\geq1:\; \sigma(2n+1)\geq \sigma(2n)\}$ and let $B(x) = |\{n\in B: \;  n\leq x\}|$. Then $\dens B$ exists and we have
\begin{equation}\label{rose}
\lowerbound \leq \dens B \leq \upperbound.
\end{equation}
\end{thm}

The precision in (\ref{rose}) is not as high as in the analogous problem concerning abundant numbers, that is, those numbers $n$ such that $\sigma(n)/n \geq 2$. Let $\dens A$ be the natural density of abundant numbers. We have that $0.247617 < \dens A < 0.247648$, due to the first author \cite{MK1,MKThesis}. We shall draw on methods used in \cite{Deleglise, Mattics} to establish Theorem \ref{flower}.

In \S \ref{stomach} we prove that the density of $B$ exists. In \S \ref{intestine} we set up the tools to bound $\dens B$ and in \S \ref{foot} we complete the proof of Theorem \ref{flower}.

\section{Existence of $\dens B$} \label{stomach}
Let $h(n)=\sigma(n)/n$. It will be convenient to work with the set
\[C:=\{n:h(2n+1)\geq h(2n)\}.\]
We will prove that the sets $B$ and $C$ have equal densities. 
First observe that
\[\frac{h(2n+1)}{h(2n)}=\frac{\sigma(2n+1)}{\sigma(2n)}\cdot\frac{2n}{2n+1},\]
so $C\subseteq B$. By $\cite{Shapiro}$, $C$ has a density, so it remains to prove that the set
\[ B-C=\left\{n : 0\leq \sigma(2n+1)-\sigma(2n)<\frac{\sigma(2n)}{2n}\right\} \]
has density zero. On the one hand, Gr\"onwall's theorem \cite{Gronwall} states that
\[\limsup_{n\to\infty}\frac{\sigma(n)/n}{\log\log n}=e^\gamma, \]
where $\gamma$ is the Euler--Mascheroni constant. Hence, for $n\in B-C$ we have that
\[ \sigma(2n+1)-\sigma(2n)=O(\log\log n). \]
On the other hand, Lemma 2.1 of \cite{LP} gives that on a set $S$ of asymptotic density 1, $p\mid \sigma(n)$ for every prime $p\leq \log\log n/\log\log\log n$. Writing $K(n)$ for the product of the primes satisfying this inequality, the prime number theorem yields
\[ K(n)=\log(n)^{(1+o(1))/\log\log\log n}. \]
Thus, for almost all $n$, $K(2n)\mid \sigma(2n+1)-\sigma(2n)$. It follows that in set $B-C$,
either $\sigma(2n+1)=\sigma(2n)$ or
\[ \log(n)^{(1+o(1))/\log\log\log n}=K(2n)\leq \sigma(2n+1)-\sigma(2n)=O(\log\log n),  \]
a contradiction for sufficiently large $n$. In the case of equality, we invoke the result in \cite{Erdos} or \cite{EPS} that the set of $n$ satisfying the equality has density zero. This establishes that the set $B$ has a density and that $\dens B=\dens C$.

\section{Preparatory results}\label{intestine}
In this section, we partition the set $C$ into subsets and bound the densities of these subsets. 

\subsection{Smooth partitions}
Let $y\geq 2$. We say a number $n$ is $y$-smooth if its largest prime divisor $p$ has $p\leq y$, and write $S(y)$ for the set of $y$-smooth numbers. Let $Y(n)$ be the largest $y$-smooth divisor of $n$. We define
\[S(a,b):=\{n\in\bbN: Y(2n+1)=a, Y(2n)=b\}.\]
Note that the sets $S(a,b), a,b\in S(y)$ partition $\bbN$, and that $S(a,b)=\varnothing$ unless $b$ is even and $\gcd(a,b)=1$. We partition $C$ via $C(a,b):=C\cap S(a,b)$.

We will express bounds of $\dens C(a,b)$ in terms of $\dens S(a,b)$. To see that $S(a,b)$ has a natural density and to determine the value of the density, we will show that $S(a,b)$ is a finite union of arithmetic progressions. Denote the set of totatives modulo $N$ by
\[\Phi(N):=\{t\in\bbN : 1\leq t\leq N, \gcd(t,N)=1\}.\]
We define $P=P(y)$ as the product of primes $p$, $p\leq y$. For any $n\in\bbN$ we have $\gcd(n/Y(n),P)=1$, so we may partition $S(a,b)$ by
\[ S(a,b\,;t_1,t_2):=\{n\in S(a,b) : (2n+1)/a\equiv t_1 \bmod P,\, 2n/b\equiv t_2 \bmod P\},\]
for $t_1, t_2\in \Phi(P)$. We will show that these sets are either empty or are arithmetic progressions.

For $n\in S(a,b\,;t_1,t_2)$, the condition $n\in S(a,b)$ implies $2n+1=ax, 2n=by$ for some $x,y\in\bbZ$. We thus study the linear Diophantine equation
\begin{equation}\label{firstde}
    ax-by=1.
\end{equation}
Writing the congruence conditions as
\[x=t_1+x'P,\quad y=t_2+y'P, \quad x', y'\in\bbZ, \]
the equation in \eqref{firstde} becomes
\begin{equation}\label{secondde}
    aPx'-bPy'=1-at_1+bt_2.
\end{equation}
This equation has solutions if and only if $P\mid 1-at_1+bt_2$. In this case, write $P\ell=1-at_1+bt_2$. Then \eqref{secondde} simplifies to
\begin{equation*}\label{thirdde}
    ax'-by'=\ell,
\end{equation*}
which has the general solution $x'=x_0\ell+kb, y'=y_0\ell+ka, k\in\bbZ$, where $x=x_0, y=y_0$ is a particular solution for \eqref{firstde}. We conclude that $n\in S(a,b\, ;t_1,t_2)$ has the form
\begin{align*}
    2n+1&=a(t_1+P\ell)x_0+abPk,\\
    2n&=b(t_2+P\ell)y_0+abPk,
\end{align*}
%$k\in\bbZ$, 
and any choice of $k$ such that $n\in\bbN$ puts $n$ in $S(a,b\, ; t_1, t_2)$. Thus $S(a,b\, ; t_1, t_2)$ is an arithmetic progression when nonempty and
\[\dens S(a,b\, ; t_1, t_2)=\begin{cases}0 & P\nmid 1-at_1+bt_2,\\ \frac{2}{abP} & P\mid 1-at_1+bt_2.\end{cases} \]

To determine $\dens S(a,b)$, we must count the number of ordered pairs $(t_1,t_2)$ satisfying $P\mid 1-at_1+bt_2$. We check for valid pairs modulo each prime $p\mid P$. If $p\mid a$, then $p\mid 1-at_1+bt_2$ if and only if $t_2\equiv -b^{-1} \bmod p$, so $t_1$ is free and $t_2$ is completely determined modulo $p$. Thus, there are $p-1$ ordered pairs modulo $p$. Similarly, if $p\mid b$, $p\mid 1-at_1+bt_2$ if and only if $t_1\equiv a^{-1} \bmod p$, so again there are $p-1$ ordered pairs modulo $p$. Finally, if $p\nmid ab$,
$p\mid 1-at_1+bt_2$ if and only if $t_2\equiv b^{-1}(at_1-1)\bmod p$. For each $t_1\in\Phi(P)$, we fail to get a valid $t_2\in\Phi(P)$ only if  $t_1 \equiv a^{-1} \bmod p$. Thus, there are $p-2$ valid ordered pairs modulo $p$. We conclude by the Chinese remainder theorem that
\[\# \{(t_1,t_2)\in \Phi(P)^2 : P\mid 1-at_1+bt_2\}=\prod_{p\mid ab}(p-1)\prod_{p\mid P, p\nmid ab}(p-2),\]
so that
\[\dens S(a,b)=\frac{2}{ab}\prod_{p\mid ab}\left(1-\frac{1}{p}\right)\prod_{\substack{p\mid P\\ p\nmid ab}}\left(1-\frac{2}{p}\right). \]

\subsection{Moments of $h(2n+1)$ and $h(2n)$}

To bound $\dens C(a,b)$ we will also need bounds on the following moments of $h(2n+1)$ and $h(2n)$ over $n\in S(a,b)$
\[ \sum_{\substack{n\leq x \\ n\in S(a,b)}}h^r(2n+1)\quad \text{and}\quad\sum_{\substack{n\leq x \\ n\in S(a,b)}}h^r(2n).\]
To this end, we prove a higher-moments analogue of the lemma in \cite{Mattics} using ideas in \cite{Deleglise}.
\begin{lem}\label{arronax}
Let
\[g(n):=\left(\frac{\sigma(n)}{n}\right)^r, \quad \rho(p^\alpha):=g(p^\alpha)-g(p^{\alpha-1}),\]
and
\[\Lambda_k(r):=\prod_{p\nmid k}\left(1+\frac{\rho(p)}{p}+\frac{\rho(p^2)}{p^2}+\cdots\right). \]
If $h$ and $k$ are given coprime positive integers, $r\geq 1$ and $x\geq 2$, then
\[ \sum_{\substack{n\leq x\\ n\equiv h\bmod k}}\left(\frac{\sigma(n)}{n}\right)^r = x\frac{\Lambda_k(r)}{k}+O((\log k)^r). \]

\end{lem}
(Note that, although we are borrowing the notation of \deleglise, our meaning for $k$ differs from his.)
\begin{proof}
We generalize the lemma in \cite{Mattics} which proves the case $r=1$. Fix a real number $r\geq 1$.  By M\"obius inversion, we express $g(n)$ as the divisor sum
\[ g(n)=\sum_{d\mid n}\rho(d), \]
where
\[ \rho(n)=\sum_{d\mid n}\mu\left(\frac{n}{d}\right)g(d). \]
Since $g$ is multiplicative, so is $\rho$, and on prime powers $p^\alpha$ we have 
\[ \rho(p^\alpha)=g(p^\alpha)-g(p^{\alpha-1}). \]
Note that $\rho$ is always positive.

If $\chi$ is a character modulo $k$, we have
\begin{align*}
\sum_{n\leq x}\chi(n)g(n)&=\sum_{n\leq x}\chi(n)\sum_{d\mid n}\rho(d)\\
&=\sum_{d\leq x}\chi(d)\rho(d)\sum_{m\leq x/d}\chi(m).
\end{align*}
If $\chi$ is non-principal, we have
\[ \sum_{n\leq x}\chi(n)g(n)=O\left(\sum_{d\leq x}\rho(d)\right). \]
If $\chi$ is the principal character, and letting a dash on a summation denote sums restricted to integers relatively prime to $k$, we have
\begin{align*}
\sum_{n\leq x}\chi(n)g(n)&=\sum_{d\leq x}{}'\rho(d)\left(\frac{\varphi(k)}{k}\frac{x}{d}+O(1)\right)\\
&=\frac{\varphi(k)}{k} x \sum_{d\leq x}{}'\frac{\rho(d)}{d}+O\left(\sum_{d\leq x}\rho(d)\right)\\
&=\frac{\varphi(k)}{k} x \sum_{d=1}^\infty{}'\frac{\rho(d)}{d}+O\left(\sum_{d>x}\frac{\rho(d)}{d}+\sum_{d\leq x}\rho(d)\right)\\
&=\frac{\varphi(k)}{k} x \Lambda_k(r)+O\left(\sum_{d>x}\frac{\rho(d)}{d}+\sum_{d\leq x}\rho(d)\right),
\end{align*}
where
\begin{equation*}\label{ned}\Lambda_k(r):=\sum_{d=1}^\infty{}'\frac{\rho(d)}{d}.
\end{equation*}
Again by multiplicativity of $\rho$, we have
\begin{equation*}\label{land}
\Lambda_k(r)=\prod_{p\nmid k}\left(1+\frac{\rho(p)}{p}+\frac{\rho(p^2)}{p^2}+\cdots\right).
\end{equation*}

Multiplying by $\overline \chi(h)$ and summing over the characters $\chi$ modulo $k$, we obtain
\begin{align*}
\sum_{\substack{n\leq x\\n\equiv h\bmod k}}g(n)&=\frac{1}{\varphi(k)}\sum_{\chi}\overline \chi(h)\sum_{n\leq x}\chi(n)g(n)\\
&=x\frac{\Lambda_k(r)}{k}+O\left(\sum_{d>x}\frac{\rho(d)}{d}+\sum_{d\leq x}\rho(d)\right).
\end{align*}

It remains to estimate the error.  Since $\sum_d\frac{\rho(d)}{d}$ is a convergent series, its tail is $o(1)$.  We now estimate
\[ \sum_{d\leq x}\rho(d). \]
We have
\begin{align*}
\sum_{d\leq x}\rho(d)&\leq \prod_{p\leq x}\left(1+\rho(p)+\rho(p^2)+\cdots\right)\\
&=\prod_{p\leq x}\lim_{\alpha\to\infty}g(p^\alpha)\\
&=\prod_{p\leq x}\left(1+\frac{1}{p-1}\right)^r\\
&=\exp\log\left(\prod_{p\leq x}\left(1+\frac{1}{p-1}\right)^r\right)\\
&=\exp\left(r\sum_{p\leq x}\log\left(1+\frac{1}{p-1}\right)\right)\\
&\leq \exp\left(r\sum_{p\leq x}\frac{1}{p-1}\right),
\end{align*}
where we have used the bound $\log(1+x)\leq x$ for $x>0$. Since
\[ \frac{1}{p-1}=\frac{1}{p}+O\left(\frac{1}{p^2}\right), \]
and
\[ \sum_{p\leq x}\frac{1}{p}=\log\log x +O(1),\quad \sum_{p\leq x}\frac{1}{p^2}=O(1), \]
we conclude that
\[\sum_{d\leq x}\rho(d)=O\left((\log x)^r\right). \]
\end{proof}

By Lemma \ref{arronax} and our characterization of the set $S(a,b\,;t_1,t_2)$ as an arithmetic progression when $P\mid 1-at_1+bt_2$, we conclude that for such pairs $(t_1,t_2)$ we have
\begin{align*}
     \sum_{\substack{n\leq x \\ n\in S(a,b\,;t_1,t_2)}}h^r(2n+1)&=h^r(a)\sum_{\substack{m\leq (2x+1)/a\\m\equiv (t_1+P\ell)x_0 \bmod bP}}h^r(m)\\
     &=h^r(a)\Lambda_P(r)\frac{2}{abP}x+O(\log^rx).
\end{align*}
Summing over all pairs $(t_1,t_2)$, we have
\[ \sum_{\substack{n\leq x \\ n\in S(a,b)}}h^r(2n+1)\sim h^r(a)\Lambda_P(r)\dens S(a,b)x,\quad x\to \infty.\]
Likewise,
\[ \sum_{\substack{n\leq x \\ n\in S(a,b)}}h^r(2n)\sim h^r(b)\Lambda_P(r)\dens S(a,b)x,\quad x\to\infty.\]

\section{Bounds on $\dens B$}\label{foot}

We can now place bounds on $\dens C$, and thus on $\dens B$, by bounding $\dens C(a,b)$. We call $0$ and $\dens S(a,b)$ 
trivial bounds for $\dens C(a,b)$. For a nontrivial upper bound, we observe that
\begin{align*}
    \sum_{\substack{n\in S(a,b)\\n\leq x}}h^r(2n+1)&=\sum_{\substack{n\in S(a,b)\\ n\in C \\n\leq x}}h^r(2n+1)+\sum_{\substack{n\in S(a,b)\\ n\not\in C \\n\leq x}}h^r(2n+1)\\
    &\geq \sum_{\substack{n\in S(a,b)\\ n\in C \\n\leq x}}h^r(2n)+\sum_{\substack{n\in S(a,b)\\ n\not\in C \\n\leq x}}h^r(2n+1)\\
    &\geq h^r(b)|C(a,b)\cap[1,x]|+h^r(a)(|S(a,b)\cap[1,x]|-|C(a,b)\cap[1,x]|).
\end{align*}
Dividing by $x$ and taking $x\to\infty$ we have
\[h^r(a)\Lambda_P(r)\dens S(a,b) \geq h^r(b)\dens C(a,b) + h^r(a)\dens S(a,b) - h^r(a)\dens C(a,b). \]
In the case $h(b)>h(a),$ we arrive at the upper bound
\[ \dens C(a,b)\leq \frac{h^r(a)(\Lambda_P(r)-1)}{h^r(b)-h^r(a)}\dens S(a,b).\]
For this upper bound to be nontrivial, we require $\Lambda_P(r)^{1/r}<h(b)/h(a).$ Note that since $\Lambda_P(r)>1$ for all $r\geq 1$, this condition implies $h(b)>h(a)$.

For a nontrivial lower bound, we proceed similarly:
\begin{align*}
\sum_{\substack{n\in S(a,b)\\n\leq x}}h^r(2n)&=\sum_{\substack{n\in S(a,b)\\ n\in C \\n\leq x}}h^r(2n)+\sum_{\substack{n\in S(a,b)\\ n\not\in C \\n\leq x}}h^r(2n)\\
&\geq \sum_{\substack{n\in S(a,b)\\ n\in C \\n\leq x}}h^r(2n)+\sum_{\substack{n\in S(a,b)\\ n\not\in C \\n\leq x}}h^r(2n+1)\\
&\geq h^r(b)|C(a,b)\cap[1,x]|+h^r(a)(|S(a,b)\cap[1,x]|-|C(a,b)\cap[1,x]|).
\end{align*}

Thus, asymptotically we have
\[h^r(b)\Lambda_P(r)\dens S(a,b) \geq h^r(b)\dens C(a,b) + h^r(a)\dens S(a,b) - h^r(a)\dens C(a,b). \]
In the case $h(b)-h(a)<0,$ we have
\[ \dens C(a,b) \geq \frac{h^r(a)- h^r(b)\Lambda_P(r)}{h^r(a)-h^r(b)}\dens S(a,b).\]
This bound is nontrivial when $h(a)/h(b)>\Lambda_P(r)^{1/r}$, and this condition implies $h(b)<h(a)$. 

For upper bounds $\Lambda_P^+(r)$ for $\Lambda_P(r)$ we use the work of Del\'eglise \cite{Deleglise} when $r>1$, where we have taken 65536 to be the maximum prime bound:
\[\Lambda_P^+(r)=\prod_{\substack{p \text{ prime}\\y<p< 65536}}\left(1+\frac{(1+1/p)^r-1}{p}+\frac{r}{(p^4-p^2)\left(1-\frac{1}{p}\right)^{r-1}}\right)\exp(1.6623114\times 10^{-6}r).
\]
When $r=1$ we use
\[\Lambda_P(1)=\Lambda_P^+(1)=\zeta(2)\prod_{p\mid P}\left(1-\frac{1}{p^2}\right). \]

To summarize, we use the following bounds for $\dens C(a,b)$:
\begin{align*}
    \dens C(a,b) & \geq \dens C^-(a,b)=\begin{cases} \frac{h^r(a)- h^r(b)\Lambda_P^+(r)}{h^r(a)-h^r(b)}\dens S(a,b) & \text{ for }h(a)/h(b)>\Lambda_P^+(r)^{1/r},\\
    0 & \text{ for }h(a)/h(b)\leq \Lambda_P^+(r)^{1/r},\end{cases}\\
    \dens C(a,b) & \leq \dens C^+(a,b)=\begin{cases} \frac{h^r(a)(\Lambda_P^+(r)-1)}{h^r(b)-h(a)}\dens S(a,b) & \text{ for }h(b)/h(a)>\Lambda_P^+(r)^{1/r},\\
    \dens S(a,b) & \text{ for }h(b)/h(a)\leq\Lambda_P^+(r)^{1/r}.\end{cases}\\    
\end{align*}

Then
\[ \sum_{a,b\in S(y)}\dens C^-(a,b)\leq \dens C \leq \sum_{a,b\in S(y)}\dens C^+(a,b).\]

In practice, we fix the parameters $y, z,$ and $r_{\text{max}}$, then recursively run through odd $a\in S(y)\cap[1,z]$. For each $a$ we recursively run through even $b\in S(y)\cap[1,z/a]$. For a given pair $(a,b)$, we calculate $C^\pm(a,b)$ for $1\leq r\leq \min(r_1,r_{\text{max}})$ where $r_1$ is the value of $r$ that produces a locally optimum bound. For example, in the case of $\dens C^+(a,b)$, we calculate bounds consecutively from $r=1$ until the values stop decreasing or we reach $r=r_{\text{max}}$, then keep the minimum value found.

By experimentation, we find that different values of the parameters $y$ and $z$ optimize the upper and lower bounds over a comparable time period. For the lower bound, the choice $y=353, z=10^{13}, r_{\text{max}}=2000$ yielded the value $0.0539171$ in $34.4$ hours. For the upper bound, the choice $y=157, z=10^{16}, r_{\text{max}}=2000$ yielded the value $0.0549446$ in $25.1$ hours. Both of these calculations were done on a Dell XPS 13 9370 laptop. This proves Theorem \ref{flower}.

\section*{Acknowledgements}
We thank Carl Pomerance for introducing the authors to each other, and the first author acknowledges the generous time that Carl spent as a sounding board and valuable resource throughout this project.

\end{document}